\documentclass[10pt,a4paper]{article}
\usepackage[latin1]{inputenc}
\usepackage[english]{babel}
\usepackage{amsfonts}
\usepackage{amssymb}
\usepackage{amsmath}
\usepackage{latexsym}
\usepackage[T1]{fontenc}
\newtheorem{theoreme}{Theorem}[section]
\newtheorem{lemme}[theoreme]{Lemma}
\newtheorem{corollaire}[theoreme]{Corollary}
\newtheorem{proposition}[theoreme]{Proposition}
\newtheorem{remarque}[theoreme]{Remark}
\newtheorem{conjecture}[theoreme]{Conjecture}
\newtheorem{defi}[theoreme]{Definition}
\newenvironment{preuve}{\emph{Proof} : }{\begin{flushright}$\Box$\end{flushright}}
\newcommand{\F}{\mathbb{F}}

\newcommand{\pgcd}{\mathrm{gcd}}

\begin{document}
\title{A proof of two conjectures on APN functions.}
\author{Elodie Leducq\footnote{IMJ, PhD thesis under direction of Jean-François Mestre}}
\date{}
\maketitle
\section{Introduction}

In \cite{Hell_new_families}, after a computer search, the authors give a complete list of power mappings almost perfectly nonlinear (APN) on $\F_{p^n}$, for $$p^n\in\{2^2,2^3,\ldots,2^{11},3^2,\ldots,3^7,5^2,\ldots,5^5,7^2,7^3,7^4,11^2,11^3\}.$$ Their article contains many theorems showing that some of these functions are members of an infinite family of APN functions. In \cite{MR2679075}, Zha and Wang give theorems that explain several new cases. In \cite{DOB_apn_odd_char}, Dobbertin, Mills, Müller, Pott et Willems try to find families containing values not yet explained in \cite{Hell_new_families}. They make conjectures that we prove here.  
\\\\First we recall the definition of APN function:

\begin{defi}Let $q=p^n$, $p$ being a prime number and $n$ an integer. Let $f:\F_{q}\rightarrow\F_{q}$; for all $a$ and $b\in\F_q$, we denote by $N_f(a,b)$ the number of solutions in $\F_q$ of the equation $f(x+a)-f(x)=b$. We say that $f$ is APN if $$\Delta_f:=\max(N_f(a,b),a,b\in\F_{q},a\neq0)=2.$$\end{defi}

\begin{remarque}

\begin{itemize}\item On $\F_{p^n}$, if we denote $\Delta_{x^d}$ by $\Delta_d$, we have $\Delta_{dp^i}=\Delta_d$, for $0\leq i\leq n-1$.
\item We also have $\Delta_d=\max(N_{x^d}(1,b),b\in\F_{q})$.\end{itemize}\end{remarque}

In the following table, we give all cases not yet explained in \cite{Hell_new_families} :
\begin{center}\begin{tabular}{|c|c|c|c|}
\hline
 & $p^n$ & $d$ &$dp^i$\\
 \hline
 I & $3^5$ & 134 & (134,160,238,230,206) \\
 II & $3^5$ & 152 & (152,214,158,232,212)\\
 III & $3^7$ & 40 & (40,120,360,1080,1054,976,742)\\
 IV & $3^7$ & 224 & (224,672,2016,1676,656,1968,1532)\\
 V & $3^7$ & 274 & (274,822,280,840,334,1002,820)\\
 VI & $5^3$ & 14 & (14,70,102)\\
 VII & $5^5$ & 843 & (843,1091,2331,2283,2043)\\
 \hline
 \end{tabular}\end{center}

We put $m=\frac{n+1}{2}$.

In \cite{DOB_apn_odd_char}, Dobbertin, Mills, Müller, Pott and Willems make conjectures which explain cases I, II, III, V and VII of the table :

\begin{conjecture}\label{2.1}For $n\geq5$ an odd integer, the function $x\mapsto x^d$ is APN over $\F_{3^n}$ for $$d=\left\{\begin{array}{ll} \frac{3^m-1}{2} & \textrm{if $n\equiv3 \mod4$}\\ \frac{3^m-1}{2}+\frac{3^n-1}{2}&\textrm{if $n\equiv1\mod4$}\end{array}\right..$$\end{conjecture}

\begin{conjecture}\label{2.2}For $n\geq5$ an odd integer, the function $x\mapsto x^d$ is APN over $\F_{3^n}$ for $$d=\left\{\begin{array}{ll} \frac{3^{n+1}-1}{8} & \textrm{si $n\equiv3 \mod4$}\\ \frac{3^{n+1}-1}{8}+\frac{3^n-1}{2}&\textrm{si $n\equiv1\mod4$}\end{array}\right..$$\end{conjecture}

\begin{conjecture}\label{5}Let $n$ be an odd integer. The function $x\mapsto x^d$ is APN over $\F_{5^n}$ for $$d=\frac{5^n-1}{4}+\frac{5^m-1}{2}.$$ \end{conjecture}

Now we recall two theorems proved in \cite{MR2679075} (theorem 4.1 and 4.4) that explain cases IV and VII of the table : 

\begin{theoreme}(Zha, Wang) On $\F_{3^n}$, the function $f:x\mapsto x^d$ satisfies $\Delta_d\leq2$ for $d$ such that $(3^k+1)d-2=u(3^n-1)$ where $u$ is odd and $\gcd(n,k)=1$. Furthermore, $f$ is APN if $2k<n$.\end{theoreme}

\begin{theoreme}\label{4.4}(Zha, Wang) The function $x\mapsto x^d$ is APN over $\F_{5^n}$ for $d$ such that $(5^k+1)d-2=u(5^n-1)$ where $\gcd(n,k)=1$, $u$ is odd and $k$ is even.\end{theoreme}

Using theorem \ref{4.4} Zha and Wang prove conjecture \ref{5}.

The following theorems are proved in \cite{DOB_apn_odd_char} (Theorem 2.1 and 2.2) :
\begin{theoreme}\label{3}(Dobbertin, Mills, Müller, Pott, Willems) Let $n$ be an odd integer. In $\F_{3^n}$, the function $x\mapsto x^d$ satisfies $\Delta_d\leq2$ for $$d=\left\{\begin{array}{ll} \frac{3^m-1}{2} & \textrm{if $n\equiv3 \mod4$}\\ \frac{3^m-1}{2}+\frac{3^n-1}{2}&\textrm{if $n\equiv1\mod4$}\end{array}\right..$$\end{theoreme}

\begin{theoreme}\label{4}(Dobbertin, Mills, Müller, Pott, Willems) Let $n$ be an odd integer. In $\F_{3^n}$, the function $x\mapsto x^d$ satisfies $\Delta_d\leq2$ for $$d=\left\{\begin{array}{ll} \frac{3^{n+1}-1}{8} & \textrm{si $n\equiv3 \mod4$}\\ \frac{3^{n+1}-1}{8}+\frac{3^n-1}{2}&\textrm{si $n\equiv1\mod4$}\end{array}\right.$$\end{theoreme}

~\\In the next section, we prove conjectures \ref{2.1} and \ref{2.2}. In section 3, we prove the following theorem, which gives a new infinite family of APN functions :

\begin{theoreme}\label{n5}Let $l\geq 2$ and $n$ an integer such that $n\equiv -1\mod2^l$ then the function $x\mapsto x^d$ is APN over $\F_{5^n}$ for $$d=\frac{1}{2}\frac{5^{n+1}-1}{5^{\frac{n+1}{2^l}}+1}+\frac{5^n-1}{4}.$$ \end{theoreme}
Finally, in the last section, we make some remarks about Zha and Wang theorems. In particular we give a new proof of conjecture \ref{5} in the case where $n\equiv3\mod4$.

\section{Proof of conjectures \ref{2.1} and \ref{2.2}}

In this part $p=3$. Using theorems \ref{3} and \ref{4} , we only have to show that $\Delta_d\neq1$ which means that $(x+1)^d-x^d$ is not a permutation polynomial. We do that with Dickson and Hermite's criterion which gives a necessary and sufficient condition for a polynomial with coefficients in $\F_q$ to be a permutation polynomial over $\F_q$ :
\begin{theoreme}(Hermite an Dickson's criterion, see \cite{lidl_finite_fields} p. 349) The following propositions are equivalent :
\begin{enumerate}\item $f(x)\in\F_q[x]$ is a permutation polynomial over $\F_q$.
\item $f(x)$ has exactly one root in $\F_q$ and
\\ $\forall t$, $1\leq t\leq q-2$, $t\not\equiv 0\mod p$, the reduction $(f(x))^t\mod x^q-x$ has degree less than $q-1$\end{enumerate}\end{theoreme}

We also need Lucas' theorem :
\begin{theoreme}(Lucas see \cite{Lucas}, p.230) Let $p$ be a prime number, $n$ and $r$ integers. We consider p-adic decomposition of $n$ and $r$ :
\\$n=n_0+n_1p+\ldots+n_kp^k$ with $0\leq n_i\leq p-1$,
\\$r=r_0+r_1p+\ldots+r_kp^k$ with $0\leq r_i\leq p-1$.
\\Then $$\binom{n}{r}\equiv\prod_{i=0}^k\binom{n_i}{r_i}\mod p.$$\end{theoreme}

Now we can prove the conjectures \ref{2.1} and \ref{2.2}.

Since in each case $\gcd(d,3^n-1)=2$, we have to find $t$ in Hermite and Dickson's criterion such that the degree of $((x+1)^d-x^d)^t\mod x^q-x$ is $q-1$. 

\begin{align*}\left((x+1)^d-x^d\right)^t&=\sum_{k=0}^t\binom{t}{k}(x+1)^{dk}(-1)^{t-k}x^{d(t-k)} \\&=\sum_{k=0}^t\binom{t}{k}(-1)^{t-k}\left(\sum_{j=0}^{dk}\binom{dk}{j}x^{dk-j}\right)x^{d(t-k)}\\&=\sum_{k=0}^t\binom{t}{k}(-1)^{t-k}\sum_{j=0}^{dk}\binom{dk}{j}x^{dt-j} \\&=\sum_{j=0}^{dt}x^{dt-j}\sum_{k=\lceil\frac{j}{d}\rceil}^t(-1)^{t-k}\binom{t}{k}\binom{dk}{j}\end{align*}
The degree of $x^k \mod x^q-x$ is $q-1$ if and only if $k\equiv0\mod q-1$ and $k\neq0$. For $v\in\mathbb{R}$ we denote by $\lceil v\rceil$ the ceilling of $x$ and by $\lfloor v\rfloor$ its floor. Then the coefficient of $x^{q-1}$ in $\left((x+1)^d-x^d\right)^t \mod x^q-x$ is $$C=\sum_{i=1}^{\lfloor\frac{dt}{q-1}\rfloor}\sum_{k=\lceil\frac{dt-i(q-1)}{d}\rceil}^t(-1)^{t-k}\binom{t}{k}\binom{dk}{dt-i(q-1)}.$$

\paragraph*{Case where $n\equiv 1\mod 4$ in conjecture \ref{2.1}}
~~\\\\We have $d=\frac{3^m-1}{2}+\frac{3^n-1}{2}$. We choose $t=2$ then, since $n\ne1$, $\lfloor\frac{dt}{q-1}\rfloor=1$  and $\lceil\frac{dt-(q-1)}{d}\rceil=1$. So $$C=\binom{3^n+3^m-2}{3^m-1}+\binom{\frac{3^n-1}{2}+\frac{3^m-1}{2}}{3^m-1}$$
On one side, $3^n+3^m-2=3^n+\displaystyle\sum_{k=1}^{m-1}2\times3^k+1$ and $3^m-1=\displaystyle\sum_{k=0}^{m-1}2\times3^k$ so by Lucas' theorem, $\binom{3^n+3^m-2}{3^m-1}\equiv 0\mod 3$. \\On the other side, $\frac{3^n-1}{2}+\frac{3^m-1}{2}=\displaystyle\sum_{k=0}^{n-1}3^k+\sum_{k=0}^{m-1}3^k=\sum_{k=m}^{n-1}3^k+\sum_{k=0}^{m-1}2\times3^k$. So $\binom{\frac{3^n-1}{2}+\frac{3^m-1}{2}}{3^m-1}\equiv1\mod 3$. 

\paragraph*{Case where $n\equiv 3\mod 4$ in conjecture \ref{2.1}}
~~\\\\We have $d=\frac{3^m-1}{2}$. Let $s=2\times(3^{m-1}+1)$, we choose $t=2s+4$. For $n>3$, $0<t<q-1$ and $t\not\equiv0\mod3$. Furthermore, $\lfloor\frac{dt}{q-1}\rfloor=2$, $\lceil\frac{dt-(q-1)}{d}\rceil=s+6$ and $\lceil\frac{dt-2(q-1)}{d}\rceil=7$.
So $$C=\sum_{k=s+6}^t(-1)^k\binom{t}{k}\binom{dk}{dt-(q-1)}+\sum_{k=7}^t(-1)^k\binom{t}{k}\binom{dk}{dt-2(q-1)}.$$
We have $t=4\times3^{m-1}+8=3^{m}+3^{m-1}+2\times3+2$.
\\So by Lucas' theorem, $\binom{t}{k}\not\equiv0\mod3$ if and only if $k=a3^{m}+b3^{m-1}+c3+d$ where $a$, $b\in\{0,1\}$ and $c$, $d\in\{0,1,2\}$.
\\Then $$dk= a3^n+(a+b)\sum_{j=m+1}^{n-1}3^j+(a+b+c)3^{m}+(b+c+d)3^{m-1}+(c+d)\sum_{j=1}^{m-2}3^j+d$$
Furthermore, $dt-2(q-1)=3^{m+1}+\displaystyle\sum_{j=1}^{m-2}2\times3^j+1$ and \\$dt-(q-1)=3^n+3^{m+1}+\displaystyle\sum_{j=1}^{m-2}2\times3^j$.
\\So $\binom{dk}{dt-(q-1)}\not\equiv0\mod3$ if $a=1$ and $c+d=2$.
\\Assume that $b=1$, then $b+c+d=3$ and $5\geq a+b+c+1\geq3$. So the coefficient of $3^{m+1}$ in 3-adic decomposition of $dk$ is 0 which means that $\binom{dk}{dt-2(q-1)}\equiv0\mod3$.
Hence the only $k$ remaining in the first sum are $3^{m}+2\times3$, $3^{m}+3+1$ and $3^{m}+2$.
\\Now we consider the second sum : $\binom{dk}{dt-2(q-1)}\not\equiv0\mod3$ if   $d\in\{1,2\}$ and $c+d=2$, namely if $c=1$ and $d=1$ or if $c=0$ and $d=2$.
Since $k\geq7$, $a$ and $b$ can't be both 0.
Assume that $a=b=1$ then $b+c+d=3$ and $4\geq a+c+b+1 \geq3$. So the coefficient of $3^{m+1}$ in 3-adic decomposition of $dk$ is 0, which means that $\binom{dk}{dt-2(q-1)}\equiv0\mod3$. Hence the only $k$ remaining in this sum are $3^{m}+3+1$, $3^{m}+2$, $3^{m-1}+3+1$ and $3^{m-1}+2$.
\\Finally :
\begin{align*}C&=-\binom{t}{3^{m}+2\times3}\binom{d(3^{m}+2\times3)}{dt-(q-1)}  -\binom{t}{3^{m}+3+1}\binom{d(3^{m}+3+1)}{dt-(q-1)}
\\&\hspace{0.5cm}-\binom{t}{3^{m}+2}\binom{d(3^{m}+2)}{dt-(q-1)}
-\binom{t}{3^{m}+3+1}\binom{d(3^{m}+3+1)}{dt-2(q-1)}
\\&\hspace{0.5cm}-\binom{t}{3^{m}+2}\binom{d(3^{m}+2)}{dt-2(q-1)}
-\binom{t}{3^{m-1}+3+1}\binom{d(3^{m-1}+3+1)}{dt-2(q-1)} \\&\hspace{3cm}-\binom{t}{3^{m-1}+2}\binom{d(3^{m-1}+2)}{dt-2(q-1)}\end{align*}
Now, we have $t=3^{m}+3^{m-1}+2\times3+2$, so by Lucas' theorem, $\binom{t}{3^{m}+2\times3}\equiv1\mod3$. Moreover $d(3^m+2\times3)=3^n+\displaystyle\sum_{j=m+2}^{n-1}3^j+2\times3^{m+1}+2\times3^{m-1}+2\sum_{j=1}^{m-2}3^j$ and $dt-(q-1)=3^n+3^{m+1}+\displaystyle\sum_{j=1}^{m-2}2\times3^j$, so $\binom{d(3^{m}+2\times3)}{dt-(q-1)}\equiv2\mod3$.
\\We do the same for all binomials and we get $C\equiv1\mod3.$
\paragraph*{Case where $n\equiv 1\mod 4$ in conjecture \ref{2.2}}
~~\\\\We have $d= \frac{3^{n+1}-1}{8}+\frac{3^n-1}{2}$. We choose $t=2$ then $\lfloor\frac{dt}{q-1}\rfloor=1$  and $\lceil\frac{dt-(q-1)}{d}\rceil=1$. So $$C=\binom{\frac{3^{n+1}-1}{4}+3^n-1}{\frac{3^{n+1}-1}{4}}+\binom{\frac{3^n-1}{2}+\frac{3^{n+1}-1}{8}}{\frac{3^{n+1}-1}{4}}$$
On one side, we have $\frac{3^{n+1}-1}{4}+3^n-1=\displaystyle\sum_{k=1}^{\frac{n-1}{2}}2\times3^{2k}+3^n+1$
and $\frac{3^{n+1}-1}{4}=2\displaystyle\sum_{k=0}^{\frac{n-1}{2}}3^{2k}$.
So by Lucas' theorem, $\binom{\frac{3^{n+1}-1}{4}+3^n-1}{\frac{3^{n+1}-1}{4}}\equiv0\mod3$
\\On the other side, $\frac{3^n-1}{2}+\frac{3^{n+1}-1}{8}=\displaystyle\sum_{k=0}^{\frac{n-1}{2}}2\times 3^{2k}+\displaystyle\sum_{k=0}^{\frac{n-1}{2}-1} 3^{2k+1}$.
\\So, $\binom{\frac{3^n-1}{2}+\frac{3^{n+1}-1}{8}}{\frac{3^{n+1}-1}{4}}\equiv1\mod3$.
\paragraph*{Case where $n\equiv 3\mod 4$ in conjecture \ref{2.2}}
~~\\\\We have $d=\frac{3^{n+1}-1}{8}$. We choose $t=26$, for $n\geq 5$, $1\leq t\leq q-2$ and $t\not\equiv 0\mod 3$. Furthermore $$\frac{dt}{q-1}=\frac{13}{4}(2\frac{3^n}{3^n-1}+1)$$
For $n\geq5$, $\frac{243}{242}\geq\frac{3^n}{3^n-1}\geq1$.
So $\lfloor\frac{dt}{q-1}\rfloor=9$.
\\We have $t=2+2\times3+2\times9$ so by Lucas' theorem, $\binom{t}{k}\not\equiv0\mod3$ if and only if $k=a+3b+9c$, $a$, $b$, $c\in\{0,1,2\}$.
\\\\In addition, $d=\displaystyle\sum_{k=0}^{\frac{n-1}{2}}3^{2k}$, so $dk=a+\displaystyle\sum_{k=1}^{\frac{n-1}{2}}(a+c)3^{2k}+\sum_{k=0}^{\frac{n-1}{2}}b3^{2k+1}+c3^{n+1}$.
\\We have $dt-(q-1)=\displaystyle\sum_{k=1}^{\frac{n-1}{2}} 2\times3^{2k}+2\times3^n+2\times3^{n+1}$. 
So, if we write $j=\alpha+3\beta$, $\alpha$, $\beta\in\{0,1,2\}$, $$dt-(j+1)(q-1)=\alpha+\beta3+\displaystyle\sum_{k=1}^{\frac{n-1}{2}} 2\times3^{2k}+(2-\alpha)3^n+(2-\beta)3^{n+1}.$$
So $\binom{dk}{dt-(j+1)(q-1)}\not\equiv0\mod3$ if and only if $a\geq\alpha$, $b\geq\beta$, $b\geq2-\alpha$, $c\geq2-\beta$ and $a+c=2$.
\\Finally, \begin{align*}C&\equiv\binom{26}{24}-\binom{26}{16}+\binom{26}{24}-\binom{26}{13}+\binom{26}{16}+\binom{26}{8}+\binom{26}{16}\\&\hspace{3,5cm}+\binom{26}{24}-\binom{26}{16}+\binom{26}{8}+\binom{26}{8}\mod3\\&\equiv1\mod3\end{align*}
In all cases we have proved that $C\not\equiv0\mod3$; so by Hermite and Dickson's criterion $\Delta_d\neq1$.

\section{Proof of theorem \ref{n5}}

We give first some preliminary results : 
\begin{lemme}\label{inverse}(see \cite{DOB_apn_odd_char} p. 97) If $\pgcd(d,q-1)=1$ then $\Delta_d=\Delta_{d^{-1}}$ where $d^{-1}$ is the inverse of $d$ modulo $q-1$.\end{lemme}

\begin{proposition}\label{HE}(see \cite{Hell_new_families}, corollary 1 p.484) Let $n$ and $k$ be integers such that $\gcd(2n,k)=1$. Then $x\mapsto x^d$ is APN over $\F_{5^n}$ for $d=\frac{5^k+1}{2}$.\end{proposition}

\begin{corollaire}\label{6}Let $n$ be a integer such that $n\equiv2^l-1\mod 2^{l+1}$ then $x\mapsto x^d$ is APN over $\F_{5^n}$ for $d=\frac{5^{\frac{n+1}{2^l}}+1}{2}$.
 \end{corollaire}

\begin{preuve} If $n\equiv2^l-1\mod2^{l+1}$ then $\frac{n+1}{2^l}\equiv1\mod2$. So $\pgcd(2n,\frac{n+1}{2^l})=1$ and by proposition \ref{HE} we get the result.
\end{preuve}

Now we are able to prove theorem \ref{n5}.
\\\\First we consider the case where $n\equiv2^l-1\mod2^{l+1}$. We have
\begin{align*}\frac{5^{\frac{n+1}{2^l}}+1}{2}d&=\frac{5^{n+1}-1}{4}+\frac{5^{\frac{n+1}{2^l}}+1}{2}\frac{5^n-1}{4}\\&= 5^n+\frac{5^n-1}{4}\frac{5^{\frac{n+1}{2^l}}+3}{2}\end{align*}
but if $n\equiv2^l-1\mod2^{l+1}$, $\frac{n+1}{2^l}\equiv1\mod 2$ and $5^{\frac{n+1}{2^l}}+3\equiv 0\mod8$.
So $$\frac{5^{\frac{n+1}{2^l}}+1}{2}d\equiv1\mod 5^n-1$$
Hence $d$ is invertible modulo $5^n-1$ and $d^{-1}=\frac{5^{\frac{n+1}{2^l}}+1}{2}$. By corollary \ref{6}, $x\mapsto x^{d^{-1}}$ is APN over $\F_{5^n}$, so by lemma \ref{inverse}, we get the result for $n\equiv2^l-1\mod2^{l+1}$.
\\\\If $n\equiv-1\mod 2^{l+1}$ then $(2^l+1)n\equiv2^l-1\mod2^{l+1}$. We put $$e=\frac{1}{2}\frac{5^{(2^l+1)n+1}-1}{5^{\frac{(2^l+1)n+1}{2^l}}+1}+\frac{5^{(2^l+1)n}-1}{4}.$$ By the first case $x\mapsto x^e$ is APN over $\F_{5^{(2^l+1)n}}$. Furthermore,
\begin{align*}e&=\frac{1}{2}\sum_{k=0}^{2^l-1}(-1)^{k+1}(5^{\frac{(2^l+1)n+1}{2^l}})^k+\frac{5^{(2^l+1)n}-1}{4}\\&=\frac{1}{2}\sum_{k=0}^{2^l-1}(-1)^{k+1}5^{nk}(5^{\frac{n+1}{2^l}})^k+5^n\frac{5^{2^ln}-1}{4}+\frac{5^n-1}{4}\\&=5^n(5^n-1)\frac{1}{4}\sum_{j=0}^{2^l-1}5^{nj}+\frac{5^n-1}{4}+(5^n-1)\frac{1}{2}\sum_{k=1}^{2^l-1}(-1)^{k+1}\sum_{j=0}^{k-1}5^{nj}(5^{\frac{n+1}{2^l}})^k\\&\hspace{2cm}+\frac{1}{2}\sum_{k=0}^{2^l-1}(-1)^{k+1}(5^{\frac{n+1}{2^l}})^k\end{align*}
Since $\displaystyle\sum_{k=0}^{2^l-1}5^{nj}\equiv 0\mod 4$ and $\displaystyle\sum_{k=1}^{2^l-1}(-1)^{k+1}\sum_{j=0}^{k-1}5^{nj}(5^{\frac{n+1}{2^l}})^k\equiv0\mod 2$ for $l\geq2$, we have :
$$e\equiv\frac{1}{2}\frac{5^{n+1}-1}{5^{\frac{n+1}{2^l}}+1}+\frac{5^n-1}{4}\mod5^n-1$$
Since $x\mapsto x^e$ is APN over $\F_{5^{(2^l+1)n}}$, for all $b\in\F_{5^n}$ the equation $(x+1)^e-x^e=b$ has at most two solutions in $\F_{5^{(2^l+1)n}}\supset\F_{5^n}$. So, since $e\equiv d\mod 5^n-1$, for all $b\in\F_{5^n}$, $(x+1)^d-x^d=b$ has at most two solutions in $\F_{5^n}$ and $\Delta_d\leq 2$ for $n\equiv -1\mod 2^{l+1}$.
\\Furthermore since $n$ is odd, $\frac{5^n-1}{4}\equiv1\mod2$. In addition $\displaystyle\sum_{k=0}^{2^l-1}(-1)^{k+1}(5^{\frac{n+1}{2^l}})^k\equiv 0\mod4$, so $d\equiv1\mod2$. Then both 0 and -1 are solutions of $(x+1)^d-x^d=0$ and we get the result for $n\equiv-1\mod2^{l+1}$.

\section{Some remarks about Zha and Wang theorems}

First we notice that theorem \ref{4.4} gives another proof of the case $n\equiv-1\mod2^{l+1}$ of theorem \ref{n5}.
\\\\Now, we give another proof of conjecture \ref{5}. In the case where $n\equiv1\mod 4$ it is the same that in theorem 4.5 of Zha and Wang. But the proof of the case where $n\equiv 3\mod 4$ doesn't use theorem 4.4.
\\\\First we recall the idea of the proof in case where $n\equiv1\mod4$ :  
$$\frac{5^{m}+1}{2}d\equiv 1\mod (5^n-1).$$
So by lemma \ref{inverse} and corollary \ref{6} we get the result.
\\\\ If $n\equiv 3\mod 4$ then $3n\equiv 1\mod 4$. So
\begin{align*}\frac{5^{3n}-1}{4}+\frac{5^{\frac{3n+1}{2}}-1}{2}-\frac{5^n-1}{4}-\frac{5^{m}-1}{2}&=(5^n-1)\frac{2.5^m+5^n(5^n+1)}{4} \\&\equiv 0\mod (5^n-1)\end{align*}
So $e=\frac{5^{3n}-1}{4}+\frac{5^{\frac{3n+1}{2}}-1}{2}\equiv\frac{5^n-1}{4}+\frac{5^{m}-1}{2}\mod 5^n-1$ and, by the case where $n\equiv1\mod4$, $x\mapsto x^e$ is APN over $\F_{5^{3n}}$ and $\Delta_d\leq2$ for $n\equiv3\mod4$. Furthermore $d$ is odd, so $(0+1)^d-0^d=1=(-1+1)^d-(-1)^d$ and $\Delta_d=2$.
\\\\We finish by a remark on theorem 4.1 of \cite{MR2679075} :
\\For $l\in\mathbb{N}$ and $n\equiv-1\mod 2^l$, if we take $k=\frac{n+1}{2^l}$ and $u=3$ in theorem \ref{4.4}, we get that $x\mapsto x^d$ is APN for  $$d=\frac{3^{n+1}-1}{3^{\frac{n+1}{2^l}}+1}.$$
For $l=2$, this gives an explicit family to explain case IV of the table. Actually this family contains $656=224\times 5^4$.

\nocite{Li_perm_pol_deg_6_7}

\bibliographystyle{plain} 
\bibliography{bibliothese}
\end{document}